\documentclass[10pt]{amsart}
\usepackage{amsmath,amssymb,latexsym,soul,cite}
\usepackage{color,enumitem,graphicx, tikz, mathtools, microtype}
\usepackage{fancyvrb}
\usepackage{fancyhdr}
\usepackage[left=2.65cm,right=2.65cm,top=3cm,bottom=3cm]{geometry}

\usepackage[hyperref]{backref}
\usepackage{color,enumitem,graphicx}
\usepackage[colorlinks=true,urlcolor=green,citecolor=green,linkcolor=green,linktocpage,pdfpagelabels,bookmarksnumbered,bookmarksopen]{hyperref}
\usepackage[english]{babel}

\numberwithin{equation}{section}
\newtheorem{theorem}{Theorem}[section]
\newtheorem{lemma}[theorem]{Lemma}

\theoremstyle{definition}
\newtheorem{remark}[theorem]{Remark}

\newcommand{\eps}{\varepsilon}
\newcommand{\psl}{(-\Delta)^s_p}
\newcommand{\beq}{\begin{equation}}
\newcommand{\eeq}{\end{equation}}
\newcommand{\T}{{\rm Tail}}
\newcommand{\R}{{\mathbb R}}
\newcommand{\N}{{\mathbb N}}

\def\Xint#1{\mathchoice
{\XXint\displaystyle\textstyle{#1}}%
{\XXint\textstyle\scriptstyle{#1}}%
{\XXint\scriptstyle\scriptscriptstyle{#1}}%
{\XXint\scriptscriptstyle\scriptscriptstyle{#1}}%
\!\int}
\def\XXint#1#2#3{{\setbox0=\hbox{$#1{#2#3}{\int}$ }
\vcenter{\hbox{$#2#3$ }}\kern-.6\wd0}}

\def\dashint{\Xint-}

\title[Global regularity for $\psl$]{A note on global regularity for the weak solutions of fractional $p$-Laplacian equations}

\author[A.\ Iannizzotto]{Antonio Iannizzotto}
\author[S.\ Mosconi]{Sunra Mosconi}
\author[M.\ Squassina]{Marco Squassina}

\address[Antonio Iannizzotto]{Dipartimento di Matematica e Informatica
\newline\indent
Universit\`a degli Studi di Cagliari
\newline\indent
Viale L.\ Merello 92, 09123 Cagliari, Italy}
\email{antonio.iannizzotto@unica.it}

\address[Sunra Mosconi and Marco Squassina]{Dipartimento di Informatica
\newline\indent
Universit\`a degli Studi di Verona
\newline\indent
Strada le Grazie 15, I-37134 Verona, Italy}
\email{sunramosconi@gmail.com, marco.squassina@univr.it}

\thanks{Conferenza tenuta al XXV Convegno Nazionale di Calcolo delle Variazioni, Levico 2--6 febbraio 2015.}

\begin{document}

\begin{abstract}
We consider a boundary value problem driven by the fractional $p$-Laplacian operator with a bounded reaction term. By means of barrier arguments, we prove H\"older regularity up to the boundary for the weak solutions, both in the singular ($1<p<2$) and the degenerate ($p>2$) case.
\end{abstract}

\maketitle


\section{Introduction and main result}\label{sec1}

\noindent
In this short note we summarize some new results, whose full proofs are displayed in the forthcoming paper \cite{IMS1}. Our aim is to study regularity of the minimizers of the functional
\[u\mapsto\frac{1}{p}\int_{\R^N\times\R^N}\frac{|u(x)-u(y)|^p}{|x-y|^{N+ps}}\,dx\,dy-\int_\Omega fu\,dx\]
over the functions $u\in W^{s,p}(\R^N)$ such that $u=0$ a.e. in $\Omega^c$. Here and in the sequel, $\Omega\subset\R^N$ ($N>1$) is a bounded domain with a $C^{1,1}$ boundary $\partial\Omega$, $p\in(1,\infty)$ and $s\in(0,1)$ are real numbers, and $f\in L^\infty(\Omega)$. The notion of minimizer of the above functional corresponds to that of a weak solution of the Dirichlet-type boundary value problem
\beq\label{bvp}
\begin{cases}
\psl u=f & \text{in $\Omega$} \\
u=0 & \text{in $\Omega^c$,}
\end{cases}
\eeq
where $\psl$ is the $s$-fractional $p$-Laplacian operator, defined pointwisely for sufficiently smooth $u$'s by
\beq\label{psl}
\psl u(x)=2\lim_{\eps\to 0^+}\int_{B_\eps^c(x)}\frac{|u(x)-u(y)|^{p-2}(u(x)-u(y))}{|x-y|^{N+ps}}\,dy.
\eeq
The operator in problem \eqref{bvp} is both non-local and non-linear. It embraces, as a special case, the well-known fractional Laplacian $(-\Delta)^s$ ($p=2$), as well as non-linear singular ($p<2$) and degenerate ($p>2$) cases.
\vskip2pt
\noindent
Interior regularity results are not new for problems of the type \eqref{bvp}: see for instance \cite{BP}, \cite{DKP1}, \cite{KMS}, and \cite{L}. Most results, providing H\"older regularity even for more general operators, are based on Caccioppoli-type or logarithmic estimates, non-local Harnack inequalities, and (in the case of \cite{L}) a Krylov-type approach. Interior and boundary regularity results involving fully non-linear, uniformly elliptic non-local operators, are obtained in \cite{CS} and \cite{RS1} respectively. Boundary regularity for {\em degenerate} or {\em singular} problems such as \eqref{bvp}, on the other hand, is still a {\em terra incognita}.
\vskip2pt
\noindent
In the linear case $p=2$ with $f\in L^\infty(\Omega)$, the global behavior of solutions is well understood. In particular, we focus on the results of \cite{RS}: $C^s(\overline\Omega)$ regularity is obtained for the weak solutions, and is proved to be optimal by means of explicit examples, while higher regularity, namely $C^{\beta}(\Omega)$ for any $\beta\in(0,2s)$, is achieved in the interior. Furthermore, a detailed analysis of the boundary behavior of the solution $u$ reveals that $u/\delta^s$ is H\"older continuous as well in $\overline\Omega$, where
\[\delta(x)={\rm dist}(x,\Omega^c).\]
The study of boundary regularity is particularly important in view of applications to problems with a non-linear reaction $f(x,u)$, as it allows to prove fractional versions of the Pohozaev identity \cite{RS2} and of the Brezis-Nirenberg characterization of local minimizers in critical point theory \cite{IMS, ILPS}.
\vskip2pt
\noindent
Our long-term aim is to extend to the non-linear case $p\neq 2$ these latter results, which is also the reason why we focus on weak solutions rather than other types of generalized (e.g. viscosity) solutions. A first, but important, step towards such aim consists in proving global H\"older regularity for problem \eqref{bvp}.
\vskip2pt
\noindent
Our main result is the following:

\begin{theorem}\label{main}
There exist $\alpha\in(0,1)$ and $C_\Omega>0$, depending only on $N$, $p$, $s$, with $C_\Omega$ also depending on $\Omega$, such that, for all weak solution $u\in W^{s,p}(\R^N)$ of problem \eqref{bvp}, $u\in C^\alpha(\overline\Omega)$ and
\[\|u\|_{C^\alpha(\overline\Omega)}\le C_\Omega\|f\|_{L^\infty(\Omega)}^\frac{1}{p-1}.\]
\end{theorem}

\noindent
Our method differs from those of the aforementioned papers by the fact that we do not use 'hard' regularity theory, but we prefer to employ rather elementary methods based on barriers, a comparison principle from \cite{LL}, and a special 'non-local lemma' describing how $\psl u$ changes in the presence of a perturbation of $u$ supported away from $\Omega$. We shall divide our study in two main steps:
\begin{itemize}[leftmargin=0.7cm]
\item[{\bf (a)}] Interior regularity: we prove a weak Harnack inequality for positive solutions, then we localize it and develop a strong induction argument to achieve local $C^\alpha$ bounds ($\alpha\in (0,s)$) with a multiplicative constant which may blow up approaching $\partial\Omega$;
\item[{\bf (b)}] Boundary regularity: we find an explicit solution for $\psl u=0$ on the half-space, then, by means of a variable change, we produce an upper barrier near $\partial\Omega$ and by comparison we estimate $u$ by a multiple of $\delta^s$ near $\partial\Omega$. Note that, due to the non-linear nature of the problem, we cannot use fractional Kelvin transform.
\end{itemize}
Step {\bf (b)} allows us to stabilize the constant of step {\bf (a)} as we approach the boundary, thus yielding the conclusion. In view of possible future developments, we remark here that the non-optimal H\"older exponent $\alpha$ is the outcome of {\em interior} regularity rather than of the boundary behavior. In fact, it is reasonable to expect $C^s$-interior regularity, which would ensure $C^s$-global regularity at once.

\section{Preliminary results and notation}\label{sec2}

\noindent
Let $U\subseteq\R^N$ be open, not necessarily bounded. First, for all measurable $u:U\to\R$ we define the Gagliardo seminorm
\[
[u]_{W^{s,p}(U)}^p=\int_{U\times U}\frac{|u(x)-u(y)|^p}{|x-y|^{N+ps}}\,dx\,dy,\]
then we introduce some Sobolev-type function spaces:
\[W^{s,p}(U)=\big\{u\in L^p(U):\,[u]_{W^{s,p}(U)}<\infty\big\}, \ \|u\|_{W^{s,p}(U)}=\|u\|_{L^p(U)}+[u]_{W^{s,p}(U)},\]
\[W^{s,p}_0(U)=\big\{u\in W^{s,p}(\R^N):\,u=0 \,\, \text{a.e. in $U^c$}\big\}, \ \|u\|_{W^{s,p}_0(U)}=[u]_{W^{s,p}(U)},\]
\[W^{-s,p'}(U)=\big(W^{s,p}_0(U)\big)^*, \ \frac 1 p +\frac{1}{p'}=1,\]
\[\widetilde W^{s,p}(U)=\Big\{u\in L^p_{\rm loc}(\R^N):\, u\in W^{s,p}(V) \, \text{for some $V\Supset U$,} \, \int_{\R^N}\frac{|u(x)|^{p-1}}{(1+|x|)^{N+ps}}\,dx<\infty\Big\}.\]
If $U$ is unbounded, then the space $\widetilde W^{s,p}_{\rm loc}(U)$ contains the functions $u\in L^p_{\rm loc}(\R^N)$ such that $u\in\widetilde W^{s,p}(U')$ for all $U'\Subset U$. The non-local tail for a measurable $u:\R^N\to\R$ outside a ball $B_R(x)$ is
\[\T(u;x,R)=\Big(R^{ps}\int_{B_R^c(x)}\frac{|u(y)|^{p-1}}{|x-y|^{N+ps}}\,dy\Big)^\frac{1}{p-1}.\]
For all bounded $U$ and all $\Lambda\in W^{-s, p'}(U)$, by a {\em weak solution} of $\psl u=\Lambda$ in $U$ we will mean a function $u\in\widetilde W^{s,p}(U)$ such that, for all $\varphi\in W^{s,p}_0(U)$,
\beq\label{ws}
\int_{U\times U}\frac{(u(x)-u(y))^{p-1}(\varphi(x)-\varphi(y))}{|x-y|^{N+ps}}\,dx\,dy=\Lambda(\varphi)
\eeq
(we denote $a^{p-1}=|a|^{p-2}a$ for all $a\in\R$). We remark that the left-hand side of \eqref{ws} is finite and continuous with respect to $\varphi$, since $u\in\widetilde W^{s,p}(U)$. If $U$ is unbounded, $u\in \widetilde W^{s,p}_{\rm loc}(U)$ is a weak solution if it is so in any $U'\Subset U$. Corresponding notions of {\em weak super-} and {\em sub-solution} can be given. 
\vskip2pt
\noindent
Though we are mainly concerned with weak solutions, we also define a notion of {\em strong solution}: if $f\in L^1_{\rm loc}(U)$, $u\in\widetilde W^{s,p}(U)$ is said to solve $\psl u=f$ strongly in $U$ if
\begin{equation}
\label{strongsolution}2\int_{B_\eps^c(x)}\frac{(u(x)-u(y))^{p-1}}{|x-y|^{N+ps}}\,dy\to f \ \text{in $L^1_{\rm loc}(U)$.}
\end{equation}
Now we introduce two major tools for our results. The first lemma enlightens a consequence of the non-local character of $\psl$:

\begin{lemma}\label{nll}
{\rm (Non-local lemma)} Let $u\in\widetilde W^{s,p}_{\rm loc}(U)$ be a weak (resp. strong) solution of $\psl u=f$ in $U$, with $f\in L^1_{\rm loc}(U)$, and $v\in L^1_{\rm loc}(\R^N)$ satisfy
\[{\rm dist}\big({\rm supp}(v),U)>0, \qquad \int_{U^c}\frac{|v(x)|^{p-1}}{(1+|x|)^{N+ps}}\,dx<\infty.\]
Set for a.e. Lebesgue point $x\in U$ for $u$
\[h(x)=2\int_{{\rm supp}(v)}\frac{(u(x)-u(y)-v(y))^{p-1}-(u(x)-u(y))^{p-1}}{|x-y|^{N+ps}}\,dy.\]
Then, $u+v\in\widetilde W^{s,p}_{\rm loc}(U)$ and $\psl(u+v)=f+h$ weakly (resp. strongly) in $U$.
\end{lemma}

\noindent
Another important tool is the following, whose proof follows almost immediately from \cite{LL}:

\begin{theorem}\label{cp}
{\rm (Comparison principle)} Let $U$ be bounded, and $u,v\in\widetilde W^{s,p}(U)$ satisfy $u\le v$ a.e. in $U^c$ and
\[\int_{U\times U}\frac{(u(x)-u(y))^{p-1}(\varphi(x)-\varphi(y))}{|x-y|^{N+ps}}\,dx\,dy\le\int_{U\times U}\frac{(v(x)-v(y))^{p-1}(\varphi(x)-\varphi(y))}{|x-y|^{N+ps}}\,dx\,dy\]
for all $\varphi\in W^{s,p}_0(U)$, $\varphi\ge 0$ in $U$. Then, $u\le v$ a.e. in $U$.
\end{theorem}

\begin{remark}\label{smooth}
The pointwise definition of $\psl u$, even if $u$ is smooth, is a delicate issue in the singular case. Roughly speaking, if $p\ge 2$, whenever $u\in C^2(U)$ the limiting procedure in \eqref{strongsolution} is well defined, so that formula \eqref{psl} makes sense. If $p<2$, on the other hand, such a representation is possible only for $s<2(p-1)/p$, and in fact explicit examples can be detected, of very smooth functions $u\in C^\infty_c(\R^N)$ such that the integral in \eqref{psl} does not converge at a given point. This is a well known drawback of the viscosity solution approach for singular nonlinear equations.
\end{remark}

\section{Interior regularity}\label{sec3}

\noindent
We will address the interior regularity problem with a simple proof peculiar to non-local problems. It can be seen as an analogue for "divergence form" non-local equation of the elementary proof of \cite{S}, which is restricted to non-local operators in "non-divergence" form. We begin with a weak Harnack inequality for globally non-negative super-solutions (all balls are intended as centered at $0$, except when otherwise specified).

\begin{theorem}\label{whi}
{\em (Weak Harnack inequality)} Let $u\in\widetilde W^{s,p}(B_{R/3})$ satisfy
\[\begin{cases}
\psl u\ge -K & \text{weakly in $B_{R/3}$} \\
u\ge 0 & \text{in $\R^N$,}
\end{cases}\]
for some $R>0$. Then,
\begin{equation}
\label{temp}\inf_{B_{R/4}}u\ge\sigma\Big(\dashint_{B_R\setminus B_{R/2}} u^{p-1}\,dx\Big)^\frac{1}{p-1}-C(KR^{ps})^\frac{1}{p-1}
\end{equation}
with $\sigma\in (0,1)$, $C>0$, depending only on $N,s,p$.
\end{theorem}

\noindent
{\em Sketch of the proof.} For simplicity we consider only the case $p\geq 2$, and by scaling we can also assume $R=1$. We produce a lower barrier for $u$, as follows. Pick a cut-off $\varphi\in C^\infty(\R^N)$ taking values in $[0,1]$, such that $\varphi=1$ in $B_{1/4}$, $\varphi=0$ in $B_{1/3}^c$, and $\psl\varphi$ is bounded in $B_1$ in a weak sense. We choose $\sigma\in(0,1)$ (to be determined later) and set
\[L=\Big(\dashint_{B_1\setminus B_{1/2}} u^{p-1}\,dx\Big)^\frac{1}{p-1}, \ \ w=\sigma L\varphi+\chi_{B_1\setminus B_{1/2}}u.\]
Applying Lemma \ref{nll} and the elementary inequality
\[a^{p-1}-(a-b)^{p-1}\geq 2^{2-p}b^{p-1}, \ a, b\geq 0,\]
we get weakly in $B_{1/3}$
\begin{align*}
\psl (\sigma L\varphi)(x)-\psl w(x) &= \int_{B_1\setminus B_{1/2}}\frac{\big(\sigma L\varphi(x)\big)^{p-1}-\big(\sigma L\varphi(x)-u(y)\big)^{p-1}}{|x-y|^{N+ps}}\, dy\\
&\geq c\int_{B_1\setminus B_{1/2}}\frac{u^{p-1}(y)}{|x-y|^{N+ps}}\, dy\geq cL^{p-1}.
\end{align*}
By the homogeneity properties of $\psl$ we thus have
\[\psl w\le \|\psl \varphi\|_\infty(\sigma L)^{p-1}-cL^{p-1}=L^{p-1}\big(\sigma^{p-1}C_\varphi-c\big)\]
weakly in $B_{1/3}$. Therefore either \eqref{temp} is trivial for large $C$ due to $u\geq 0$ in $\R^N$, or for suitably small $\sigma>0$ we have
\[\begin{cases}
\psl w\le -K\le\psl u & \text{weakly in $B_{1/3}$} \\
w\le u & \text{in $B_{1/3}^c$,}
\end{cases}\]
hence by Theorem \ref{cp} $w\le u$ in $B_{1/3}$, and in particular $u\ge\sigma L$ in $B_{1/4}$, which gives \eqref{temp}. 
\qed

\begin{remark}
It is worth noting that, despite the proof being quite elementary, the constant $\sigma$ in the previous weak Harnack inequality degenerates as $s\to 1$. This is due to the fact that, in the previous proof,
\[
C_\varphi=\|\psl \varphi\|_\infty\simeq \frac{1}{1-s}
\]
as $s\to 1$, however regular $\varphi$ is. This gives $\sigma\simeq (1-s)^{1/(p-1)}$ and as a consequence, all the following H\"older estimates blow up as $s\to 0$. More involved proofs (see e.g. \cite{DKP1}), closely related to the classical regularity approach for local non-linear variational equations, can however give H\"older estimates which are stable when $s\to 1$.
\end{remark}

\noindent
Now we can prove a local H\"older estimate for bounded weak solutions on a ball:

\begin{theorem}\label{lhr}
{\em (Local H\"older regularity)} Let $u\in\widetilde W^{s,p}(B_{2R})\cap L^\infty(B_{2R})$ satisfy $|\psl u|\le K$ weakly in $B_{2R}$, $R>0$. Then
\[[u]_{C^\alpha(B_R)}\le C\Big((KR^{ps})^\frac{1}{p-1}+\|u\|_{L^\infty(B_{2R})}+\T(u;0,2R)\Big)R^{-\alpha},\]
with $\alpha\in (0,s]$ and $C>0$ depending on $N,s,p$, where $[u]_{C^{\alpha}(B_R)}$ is the $C^\alpha(B_R)$-seminorm.
\end{theorem}

\noindent
{\em Sketch of the proof.} First we need to localize Theorem \ref{whi}, that is, to prove a weak Harnack inequality for super-solutions which are non-negative in a ball only. If $v\in\widetilde W^{s,p}(B_{R/3})$ satisfies
\[\begin{cases}
\psl v\ge -K & \text{weakly in $B_{R/3}$} \\
v\ge 0 & \text{in $B_R$,}
\end{cases}\]
then we may apply Lemma \ref{nll} to $v_+$ (its positive part), producing a tail term depending on $v_-$ (the negative part). Using Theorem \ref{whi}, we see that for all $\eps>0$ there exists $C_\eps>0$ (depending also on $N,p,s$) such that
\beq\label{lhi}
\inf_{B_{R/4}}v\ge\sigma\Big(\dashint_{B_R\setminus B_{R/2}} v^{p-1}\,dx\Big)^\frac{1}{p-1}-C(KR^{ps})^\frac{1}{p-1}-\eps\sup_{B_R}v-C_\eps\T(v_-;0,R),
\eeq
again with $\sigma\in(0,1)$, $C>0$ depending only on $N,p,s$. We then use a strong induction argument to produce two sequences $(m_j)$, $(M_j)$, with $m_j$ non-decreasing and $M_j$ non-increasing, such that for all $j\in\N$
\[m_j\le\inf_{B_{R/4^j}}u\le\sup_{B_{R/4^j}}u\le M_j, \ M_j-m_j=\lambda\Big(\frac{R}{4^j}\Big)^\alpha,\]
with $\alpha\in (0,1)$ depending only on $N,s,p$ and $\lambda>0$ depending on $u$. This is done by applying \eqref{lhi} to the functions $u-m_j$, $M_j-u$ in $B_{R/4^j}$, where they are both non-negative, in the inductive step. Then, we obtain the following oscillation estimate for all $r\in(0,R)$:
\beq\label{osc}
\underset{B_r}{\rm osc}\,u\le C\Big((KR^{ps})^\frac{1}{p-1}+\|u\|_{L^\infty(B_R)}+\T(u;0,R)\Big)\frac{r^\alpha}{R^\alpha},
\eeq
with $C>0$ depending on $N,s,p$. A standard argument then provides the claimed estimate. \qed


\section{Boundary regularity and conclusion}\label{sec4}

\noindent
In this final section we turn back to weak solutions of \eqref{bvp}. First, by applying Theorem \ref{cp} to $u$ and a multiple of the weak solution $\psi\in W^{s,p}_0(B_1)$ of 
\[\begin{cases}
\psl\psi=1 & \text{in $B_1$} \\
\psi=0 & \text{in $B_1^c$},
\end{cases}\]
we prove the following:

\begin{theorem}\label{apb}
{\rm (A priori bound)} Let $u\in W^{s,p}_0(\Omega)$ satisfy $|\psl u|\le K$ weakly in $\Omega$. Then
\[\|u\|_{L^\infty(\Omega)}\le (C_dK)^\frac{1}{p-1},\]
with $C_d>0$ depending on $N,s,p$, and $d={\rm diam}(\Omega)$.
\end{theorem}

\noindent
Now we produce a local upper barrier. We set $e_N=(0,\ldots,0,1)$, $\R^N_+=\{x\in \R^N:\, x\cdot e_N>0\}$.

\begin{lemma}\label{lup}
There exists $w\in C^s(\R^N)$, $r>0$, $a\in(0,1)$, $c>1$ such that
\[\begin{cases}
\psl w\ge a & \text{weakly in $B_r(e_N)\setminus B_1$} \\
c^{-1}(|x|-1)^s\le w(x)\le c(|x|-1)^s & \text{in $\R^N$.}
\end{cases}\]
\end{lemma}

\noindent
{\em Sketch of the proof.} We divide our argument in four steps:
\begin{itemize}[leftmargin=0.7cm]
\item[{\bf (a)}] We find an explicit solution on a half-space: namely, $u_N(x)=(x_N)_+^s$ belongs in $\widetilde W^{s,p}_{\rm loc}(\R^N_+)$ and satisfies $\psl u_N=0$ both strongly and weakly in $\R^N_+$.
\item[{\bf (b)}] For all $R>0$ big enough we find a diffeomorphism $\Phi\in C^{1,1}(\R^N,\R^N)$ such that $\Phi(x)=x$ outside a ball, mapping $0$ to $x_R$ and the plane $\{x_N=0\}$ to $\partial B_R(\tilde x_R)$ locally around $0$, and satisfying further $\|\Phi-I\|_{C^{1,1}}\le C/R$.
\item[{\bf (c)}] We prove stability of $\psl$ under $C^{1,1}$ changes of variable of the type above: setting $v_N=u_N\circ\Phi^{-1}$, we see that $\psl v_N=g$ weakly in $\Phi(\R^N_+)$, for some $g\in L^\infty(\R^N)$ with $\|g\|_{L^\infty(\R^N)}<C/R$.
\item[{\bf (d)}] We truncate $v_N$ at a convenient heigth $M>0$, i.e., we set $v=\max\{v_N,M\}$, and we apply Lemma \ref{nll} to obtain $\psl v\ge b$ weakly in $\Phi(\R^N_+)\setminus B_2(x_R)$ ($b>0$).
\end{itemize}
By scaling and translating $v$, we find $w$ as required. \qed
\vskip4pt
\noindent
The barrier $w$ is used to prove an estimate of weak solutions, near $\partial\Omega$, by means of a multiple of $\delta^s$:

\begin{theorem}\label{be}
Let $u\in W^{s,p}_0(\Omega)$ satisfy $|\psl u|\le K$ weakly in $\Omega$. Then for a.e. $x\in\Omega$
\[|u(x)|\le (C_\Omega K)^\frac{1}{p-1}\delta^s(x),\]
with $C_\Omega>0$ depending on $N,s,p,$ and $\Omega$.
\end{theorem}

\noindent
{\em Sketch of the proof.} We may reduce to the case $K=1$. By Theorem \ref{apb}, the desired estimate is easily obtained away from $\partial\Omega$. Due to the reguarity of $\partial\Omega$, we can find $\rho>0$ such that in the set
\[\Omega_\rho=\{x\in\Omega:\,\delta(x)<\rho\}\]
$\delta$ decreases linearly on segments normal to $\partial\Omega$. Fix $x_1\in\Omega_\rho$ and denote $x_0$ its unique metric projection on $\partial\Omega$. Let $w$ be as in Lemma \ref{lup}. By scaling and translating $w$, we construct $\tilde w\in C^s(\R^N)$ such that $\tilde w\le C\delta^s$ on the line segment $[x_0,x_1]$, and moreover
\[\begin{cases}
\psl u\le 1\le\psl\tilde w & \text{weakly in $B_r(x_0)\cap\Omega$} \\
u\le\tilde w & \text{in $(B_r(x_0)\cap\Omega)^c$,}
\end{cases}\]
with a small $r<|x_0-x_1|$. By Theorem \ref{cp} we see that $u(x_1)\le C\delta^s(x_1)$. An analogous argument applies to $-u$, yielding the conclusion. \qed

\begin{remark}\label{delta}
As a byproduct, by arguments analogous to those displayed in steps {\bf (a)--(c)} above, we prove that, for convenient $\rho,K>0$, we have $|\psl\delta^s|\le K$ both weakly and strongly in $\Omega_\rho$. In fact, it can be proved that $\psl\delta^s\in C^\beta(\overline\Omega_\rho)$ for some $\beta\in(0,1)$, which is an interesting information, as it shows that the boundary behavior of weak solutions in the general case $p>1$ is similar to that in the linear case $p=2$ (see \cite[Lemma 3.9]{RS}).
\end{remark}

\noindent
We are now ready to conclude:
\vskip4pt
\noindent
{\em Sketch of the proof of Theorem \ref{main}.} Set $K=\|f\|_{L^\infty(\Omega)}$, so $|\psl u|\le K$ weakly in $\Omega$. By Theorem \ref{apb}, we only need to prove our estimate on the H\"older seminorm. Recalling Theorem \ref{lhr}, by a covering argument we find $\alpha\in(0,s]$ (depending only on $N,s,p$) and, for all $\Omega'\Subset\Omega$, a constant $C>0$ (depending also on $\Omega'$) such that $u\in C^\alpha(\overline\Omega')$ and $[u]_{C^\alpha(\overline\Omega')}\le CK^\frac{1}{p-1}$.
\vskip2pt
\noindent
Let $\rho>0$ be as in the proof of Theorem \ref{be}. For all $x_1\in\overline\Omega_\rho$ let $r=\delta(x_1)$. Theorem \ref{lhr} produces the following estimate:
\beq\label{bhe}
[u]_{C^\alpha(B_{r/2}(x_1))}\le C\big((Kr^{ps})^\frac{1}{p-1}+\|u\|_{L^\infty(B_r(x_1))}+\T(u;x_1,r)\big)r^{-\alpha}.
\eeq
The first term in the right-hand side of \eqref{bhe} is bounded due to $\alpha\le s$, $r\le\rho$. For the second term we invoke Theorem \ref{be} and the inequalities $\delta(x)\le 2r\le 2\rho$ for $x\in B_{r}(x_1)$ to obtain
\[|u(x)|r^{-\alpha}\le CK^\frac{1}{p-1} r^{s-\alpha}\leq CK^\frac{1}{p-1}\rho^{s-\alpha},\qquad \forall x\in B_{r}(x_1).\]
Finally, the tail term is bounded by means of Theorem \ref{be} again, together with $s$-H\"older continuity of $\delta^s$, thus we have from \eqref{bhe}
\[[u]_{C^\alpha(B_{r/2}(x_1))}\le CK^\frac{1}{p-1}\qquad \forall x_1\in \bar\Omega_\rho,\ r=\delta(x_1) \]
with $C>0$ depending on $N,s,p,$ and $\Omega$. Patching together the above estimates, we reach the conclusion. \qed

\end{document}